\input amstex
\documentstyle{amsppt}
\input amsppt.sty
\expandafter\redefine\csname logo\string @\endcsname{}

\magnification=1200

\baselineskip 20pt

\TagsOnRight

\topmatter
\NoRunningHeads
\title  On categorical approach to derived preference relations in some decision making problems
\endtitle
\author Victor~V. Rozen  \\
Department of Mathematics \\
Saratov State University, Russia \\
%Astrakhanskaja 83 46000 Saratov \\
%Russia \\
Grigori~I. Zhitomirski *\\
Department of Mathematics \\
Bar-Ilan University, Israel
%52900 Ramat-Gan \\
%Israel
\endauthor
\thanks Grigori~I. Zhitomirski, Department of Mathematics,
Bar Ilan University, Ramat Gan 52900, Israel, email:zhitomg\@macs.biu.ac.il
\endthanks

\thanks *This research is partially supported for the second of authors by THE ISRAEL SCIENCE FOUNDATION
 founded by The Israel Academy of Sciences and Humanities - Center of
 Excellence Program.
\endthanks
\abstract A structure called a decision making problem is considered. The set of outcomes (consequences) is
partially ordered according to the decision maker's preferences. The problem is how these preferences affect a
decision maker to prefer one of his strategies (or acts) to another, i.e. it is to describe so called derived
preference relations. This problem is formalized by using category theory approach and reduced to a pure
algebraical question. An effective method is suggested to build all reasonable derived preferences relations and
to compare them with each other.
\endabstract
\endtopmatter
\document
\baselineskip 20pt

Key words: Decision making problem, games with ordered outcomes, preference relations, preference functor,
monoids of relations.

\head 0. Introduction
\endhead
It is known that at first a decision making problem (DMP) was built as a structure with a real-valued payoff
function. Later on DMP appeared with a vector-valued payoff function. We consider a decision making problem with
goal structure given by a partial ordering relation (preferences of a decision maker) on the set of outcomes. At
the present moment such decision making problems have become of great interest. By the time the article had been
written and tested, the interesting paper of Dubois, Fargier, Perny and Prade (2002) appeared devoted to the same
subject and contained an extensive bibliography. It gives us an opportunity to refer a reader to this paper for
the history of the problem. We only mention such papers that concern our results immediately. In spite of the
common subjects the approaches of the present paper and the one cited above are different and results have no
intersections. Below, we describe more exactly relations between the present paper and Dubua et al.(2002)  that we
name for short {\bf QDT}.

The main problem considered in the present paper is the "transference" of preferences from the set of outcomes
(consequences) to the set of decision maker's strategies (or acts), i.e. ranking the strategies knowing that
some of consequences are better than others. Such ranking is called a derived preference relation. The point is
that such transference of notions and principles used for classical DMPs to DMPs with preferences is based
usually on intuition and leads frequently to anomalies. Further, in the next section, we explain what we mean.
The similar problem is the subject of {\bf QDT}, but there are some important differences.

The set of outcomes (consequences) is assumed in our paper to be partially ordered, while in {\bf QDT} it is
completely preordered (weakly ordered). A derived preference relation is regarded as a preorder, while in {\bf
QDT} it is assumed to be a partial order. But it is not essential because our approach works in all cases. It is
essential that the authors of {\bf QDT} formulate that "the real decision problem is to build this relation
(ranking of acts) from information regarding the likelihood of states and the decision maker preference on
consequences", while we take in account only last one. They do not explain what "to build" means but use some
axioms motivated somehow or other. On the other hand, we concentrate our attention on giving a suitable
mathematical formalization of the mentioned problem and on building all derived preference relations applying
obtained mathematical results. Thus we think that the present paper is not a rival but an addition to {\bf QDT}.

DMP considered in the paper is a structure, which components are a nonempty set $X$ of alternatives (strategies)
of a decision maker, a nonempty set of states $A$, a partially ordered set $(A, \omega )$ of outcomes and a map
$F$ that assigns an outcome $F(x,y)$ to every strategy $x$ if the state of the world is $y$. It is to mention
that the map $F$ determines a map $F\sp *$ from $X$ to the set $A\sp Y$. The set $A\sp Y$ is called the set of
acts of the decision maker. Usually, the set of acts is considered instead of the set of strategies. But we do
not suppose that every act can be realized by the considered decision maker.

A decision maker has to prefer one of his strategy to another according to his preferences in the set of
outcomes which are given by the partial order $\omega $, that is he has to "build" a relation $\rho $ on the set
of his strategies. We think that this relation has to be a preorder relation and has to be consistent with the
preference relation $\omega$ and with the realization function $F$. We shall call every such relation $\rho $ a
derivative preference relation. We look for a general principle to build the derivative preference relations for
all DMPs like for example the {\bf maxmin } principle that is the main optimality principle in the classical
decision making theory. But absence of any one-valued (functional) connection between alternatives and outcomes
makes the problem non-trivial.

Observe that a decision making problem can be considered as a game such that a player knows only his own
preferences but does not know the preferences of other players. Namely this approach is used in the monograph
Moulin (1981).

A new approach in the present paper is the use of category theory as a "methodological base" for formalizing the
problem. This approach leads to a method of constructing of derived preference relations based on algebraic
theory of monoids of binary relations. The main steps of solving the set up problem are the following ones. At
first (section 2), the conditions (axioms) are formulated for a derivative preference relation. We proceed from
minimal assumptions mentioned above and the assumption that so called Pareto-domination is a part of every
derived preference relation. From the point of view of category theory, these conditions mean the existence of
special functor (preference functor) from the category of DMPs to the category of preordered sets. In {\bf QDT},
the so called Ordinal Invariance Axiom is assumed. May be it is at the bottom of some pessimistic conclusion in
the mentioned paper. Instead of this axiom, we assume much weaker one that is fulfilled for all suggested in
literature preference relations.

Further (section 3), we show that constructing of the preference functor can be reduced to the constructing of
some special functor from the category of ordered sets to the category of preordered sets. The main results are
formulated in Theorem~1 and Theorem~2, that give an effective method to construct the required functors. Our
method leads to anomalies-free preference relations, and it turns out that practically all known ways to
specification of derivative preference relations are some concrete cases of the given construction (examples in
Section 4). Using this new method we check some properties of derived preference relations. The last section 5
shows how to use the obtained results. Particularly, all derived preference relations form a complete lattice.
In general case, some maximal and minimal preference relations are described. In all finite cases, this lattice
can be effectively constructed.

In the first section, we give all necessary definitions and some motivations for the present research.

%%%%%%%%%%%%%%%%Preliminaries
\head 1 Preliminaries and motivations
\endhead

A reader can find the necessary algebraic notions and facts in the book Gr\"{a}tzer (1979), so as some notions
and facts concerning binary relations and ordered sets. But for the convenience of a reader, we recall some
fundamental notions of relation theory.

A binary relation $\varrho $ is a set of ordered pairs, usually it is a subset of the Cartesian product $A\times
B$ of two sets, in this case it is called a binary relation between elements of these sets. The relation
$\varrho \sp {-1} = \{(b,a)\vert (a,b) \in \varrho \} $ is called the inverse relation to $\varrho $. The well
known product of binary relations $\varrho $ and $\sigma $ is denoted by $\sigma \circ \varrho$, that is $\sigma
\circ \varrho =\{(a,c)\vert \, (\exists b )( (a,b)\in \varrho \and (b,c)\in \sigma )\}$. A relation $\varrho
\subseteq A\times A $ is called reflexive if $(\forall a\in A ) (a,a) \in \varrho $, it is called transitive if
$\varrho \circ \varrho \subseteq \varrho $, that is $(a,b)\in \varrho $ and $(b,c)\in \varrho $ implies
$(a,c)\in \varrho $. A transitive and reflexive relation is called a preorder (or quasiorder) relation, it is
called a partial ordering relation if in addition it is antisymmetric, that is $(a\sb 1 ,a \sb 2 )\in \varrho
\bigcap \varrho \sp {-1} \Rightarrow a\sb 1 = a \sb 2$. Since we consider in the present paper only partial
ordering relations, further we often omit the word "partial".

For notational convenience we denote preorder and order relations by means of usual signs $\leq ,\geq $ with a
symbol of the corresponding relation on top.  The greatest element $a$ with respect to a preorder relation
$\varrho $ on a set $A$ is defined as usual, that is for all $x\in A $ we have $x\leq \sp {\varrho }a $. Of
course, if it exists it is not unique in general, but with respect to an order relation it is unique if it
exists. An element of an ordered set is called minimal or an atom if between it and the least element (if it
exists) there are no elements. Dually, a maximal or a dual atom is defined. An ordered set is called a chain if
every two elements of it are comparable ($a\sb 1 \leq a\sb 2 $ or $a\sb 2 \leq a\sb 1 $). An ordered set is
called a lattice ( complete lattice) if every its two elements have (every subset has) the least upper bound
(supremum) and the greatest lower bound (infimum).

A map $f: A\to B $ from one preorder $(A,\varrho )$  to another $(B, \sigma )$ is called isotonic if it takes
pairs from $\varrho $ to pairs from $\sigma $, that is $a\sb 1 \leq \sp {\varrho } a\sb 2 \Rightarrow f(a\sb 1 )
\leq \sp {\sigma } f(a\sb 2 ) $.

If $\lambda $ is a binary relation then its first projection, $pr\sb 1 \lambda $, is the set of elements $a$
such that there exists an element $b$ such that $(a,b)\in \lambda $, its second projection, $pr\sb 2 \lambda $,
is the first projection of the inverse relation $\lambda \sp {-1}$  and its reflexive projection, $pr\sb
{\Delta} \lambda $, is the set of all elements $a$ such that $ (a,a)\in \lambda $. The equality relation on a
set $A$ we denote by $\Delta \sb A $. Thus $pr\sb {\Delta} \lambda =pr\sb 1 (\lambda \cap \Delta \sb A)=pr\sb 2
(\lambda \cap \Delta \sb A) $ for $ \lambda \subset A\times A$.

As for the theory of categories and functors, we can refer a reader to the book Mac Lane (1971). But we use only
basic notions as the most convenient language for our aim, and it is necessary to know only basic definitions. A
category $\Cal C $ consists of two collections, $Ob(\Cal C )$, whose elements are the so called {\it objects }
of $\Cal C $,  and $Hom  (\Cal C )$, the collection of {\it morphisms } (or {\it arrows, homomorphisms, maps})
of $\Cal C $. To each morphism is assigned a pair of objects, called the {\it domain } and {\it codomain } of
the morphism. The notation $f: A\to B $ means that $f$ is a morphism with the domain $A$ and the codomain $B$.
If $f: A\to B $ and $g: B\to C $ are two morphisms, there is a morphism $g\circ f : A\to C$ called the {\it
composite } of $g$ and $f$. For every object $A$ there is a morphism $id\sb A: A\to A$, called the {\it identity
} of $A$. These data are subject to the following axioms:

\roster

(1) For $f:A\to B ,\; g: B\to C ,\; h:C\to D $ ,
$$ h\circ (g\circ f ) = (h\circ g)\circ f ;$$

(2) For $f:A\to B$,
$$ f\circ id\sb A =id\sb B \circ f =f.$$
\endroster

If $\Cal C $ and $\Cal D $ are categories, a {\it functor } $F: \Cal C \to \Cal D $ is a map for which

\roster

(1) If $f:A\to B$ is a morphism of $\Cal C $, then $F(f): F(A)\to F(B) $ is a morphism of $\Cal D $;

(2)$F(id\sb A ) = id \sb {F(A)}$;

(3) If $f:A\to B ,\; g: B\to C $ , then $F(g\circ f )=F(g) \circ F(f)$.

\endroster

All other definitions will be given when it will be necessary.

It is known that the {\bf maximin } principle is the main optimality principle in the classical decision making
theory. Let $X$ be a set of alternatives (strategies) of a decision maker, call {\bf I},  $Y$ be a set of medium
states and $f(x,y)$ be the payoff of {\bf I} if he chooses the strategy $x$ and the medium is situated in the
state $y$. The alternative $x\sb 0$ is {\bf maximin } if $\sup_{x} \inf \sb y \, f(x,y)=\inf\sb y\, f(x\sb 0
,y)$, i.e. if it maximizes the estimate $x\to \inf\sb y \, f(x,y)$.

If one considers the triplet $\Cal G = (X,Y,f)$ as a two-player zero sum game, it appears along with the {\bf
maximin } principle for the player {\bf 1} the dual one for the player {\bf 2},  i.e. the {\bf minimax }
principle: the alternative $y\sb 0$ is {\bf minimax } if it minimizes the estimate $y\to \sup\sb x\, f(x,y)$,
i.e. $\inf\sb y \, \sup\sb x\, f(x,y)=\sup\sb x\,f(x ,y\sb 0)$. The values $v\sb 1 =\sup\sb x \inf\sb y \,
f(x,y)$ and $v\sb 2 =\inf\sb y \, \sup\sb x\, f(x,y)$ are called the lower value and the upper value of the game
respectively. By the way, $v\sb 1 \leq v\sb 2$ holds always. If $v\sb 1 = v\sb 2$ holds it is said that the game
$\Cal G $ has a value.

In two-player zero sum game theory, the saddle point principle is the most important one. We recall that a point
$(x\sb 0,y\sb 0)$ is called a {\it saddle point} if $f(x,y\sb 0)\leq f(x\sb 0 ,y\sb 0)\leq f(x\sb 0 ,y)$ for all
$x\in X , y\in Y$ . The existence of a saddle point implies the existence of the game value. Conversely, if a
game has the value then every pair $(x\sb 0 ,y\sb 0)$, where $x\sb 0$ is a {\bf maximin } strategy of the player
{\bf 1 } and $y\sb 0$ is a {\bf minimax } strategy of the player {\bf 2}, forms a saddle point, and every saddle
point can be obtained by this way.

A game with ordered outcomes is a game in which preferences of players are given by partial ordering relations
on the set of outcomes, that is the goal structure of a such game is a set of outcomes that every player orders
according to his own preferences. In the case of an antagonistic game the goal structure is given by means of
one partial ordering relation that expresses the preferences of the player {\bf 1}, and the preferences of his
opponent (the player {\bf 2}) are given by the inverse partial ordering relation.

We consider an antagonistic game as an algebraic structure $\,G=(X,Y,A,\omega, F)\,$, where $X$ and $Y$ are the
sets of strategies of players {\bf 1} and {\bf 2}, $(A, \omega ) $ is a ordered set, $F$ is a map $F: X\times Y
\to A$ that is called a realization function. It is assumed only that the sets $X,Y$ and $A$ contain more than
one element. The fact $a\sb 1 \geq \sp {\omega} a\sb2 $ is interpreted in the sense that the outcome $a\sb 1$ is
more preferable for the player 1 than the outcome $a\sb 2$. If the player 1 chooses a strategy $x$ and the
player 2 chooses a strategy $y$ it comes out a situation $(x,y)$ that leads to an outcome $F(x,y)$. So, the
problem arises for every player how to prefer one of his strategy to another.

The most natural point of view is the following one: a strategy $x\sb 1$ is more preferable for the player {\bf
1} than a strategy $x\sb 2$ if  $F(x\sb 1 ,y) \geq \sp {\omega} F(x\sb2, y)$ takes place for every strategy $y$
of the player {\bf 2 }. It is so called {\it Pareto-domination }. The {\it strong Pareto-domination }  means
that all inequalities in the definition of Pareto-domination are strong. The maximal according to
Pareto-domination strategies are called {\it Pareto-optimal.} But this way has some shortcomings because the
problem is reduced to the problem of choosing one of the Pareto-optimal strategies. Nevertheless it is obvious
that if a strategy $x\sb 1$ Pareto-dominates a strategy $x\sb 2$ it is more preferable than last one in any
other kind of domination.

While searching for other more suitable principals of domination, the notions are introduced, that are analogies
to the notions of {\bf maximin } ({\bf minimax }) strategy and the value of a game in the classic case [see
Jentzsch (1964), Podinovski (1979)]. The saddle point is defined as usual (see above). Since it can be that an
{\bf infimum } of a subset of the ordered set $(A, \omega )$ does not exist, the following changing is made.

%Let $a\sp {\triangledown }$ denote the set of all elements $a'$ such that $a' \leq \sp {\omega} a $. The
%estimate $x \to {\inf} \sb y \, f(x,y) $ is replaced by the estimate
%$$
%x\to V\sb x = \bigcap_y (F(x,y))\sp {\triangledown } .
%$$
Denote by $V\sb x $ the set of all outcomes that are "guaranteed" for the player {\bf 1} by the strategy $x$,
that is: $V\sb x =\{a\in A \vert (\forall y\in Y )F(x,y)\geq\sp {\omega }a\}$. The estimate $x \to {\inf} \sb y
\, f(x,y) $ is replaced by the estimate $x\to V\sb x $. If the set $V\sb x $ is assumed as a characteristic of
the strategy $x$ the relation of so called $\alpha $-dominating for the player {\bf 1} arises:
$$ x\sb 1 \leq \sp \alpha x\sb 2\Leftrightarrow V\sb {x\sb 1} \subseteq V\sb {x\sb 2}.$$
The defined relation $\alpha $ is a preorder relation on the set $X$ of strategies of the player {\bf 1}. The
strategy $x\sb 0 \in X$ is called an $\alpha$-greatest one if it is the greatest element with respect to the
preorder relation $\alpha $. Dually, i.e. by means of rearrangement of players and replacing $\omega$ by $\omega
\sp {-1}$, the definition of an $\alpha$-greatest strategy of the player 2 is obtained.

Denote  $U\sb y =\{a \in A \,\vert (\exists x\in X )F(x,y)\geq\sp {\omega }a\}$ and set:
$$V=\bigcup_x V\sb x  \quad \text {(the lower characteristic set)};$$
$$U=\bigcap_y U\sb y  \quad \text {(the upper characteristic set)}.$$
The inclusion $V\subseteq U$ holds always. One can regard it as an analogy to the inequality $v\sb 1 \leq v\sb
2$. The equality $V=U$ is considered in Jentzsch (1964) as an analogy to the existence of a value of the game
$G$. In this case, we call the coincident lower and upper characteristic sets the {\it generalized value} of the
game $G$.

A game $\Cal G$ with payoff function can be considered formally as a game with ordered outcomes (in this case,
the set of outcomes is the set $\Bbb R$ of real numbers with respect to usual order, and the function f becomes
a payoff function). Then the $\alpha$-greatest strategy of the player 1 coincides with his {\bf maximin}
strategy, and the $\alpha$-greatest strategy of the player 2 coincides with his {\bf minimax} strategy, and the
equality $V=U$ implies the equality $v\sb 1 = v \sb 2$. Thus the introduced analogies seem to be quite natural.
But they have some anomalies (some of them were mentioned in Jentzsch (1964), Podinovski (1981), Rozen (2001)).
\medskip

{\bf Example 1}

Consider the game  $G=(X,Y,A,\omega, F)$ with ordered outcomes, where the ordered set $(A,\omega)$ is the
following 5-elements lattice: $\{0< a,b,c <1\}$, and with realization function $F$ given by the table:
$$
\vbox{ \halign{&\vrule#&
 \strut\quad\hfil#\quad\cr
&$F$&&$y_1$&&$y_2$&&$y_3$&\cr \noalign{\hrule} &$x_1$&&$b$&&$c$&&$0$&\cr \noalign{\hrule}
&$x_2$&&$0$&&$a$&&$1$&\cr} } \;.$$

Here $x\sb 1$ is the  $\alpha$-greatest strategy of the player 1, $y\sb 1 $ is the
 $\alpha$-greatest strategy of the player 2 and $V=U$. However, the pair $(x\sb 1, y\sb 1 ) $ does not form a saddle point, and this game has no saddle points at all.
Besides this, the dual condition $V\sp *=U\sp *$ is not satisfied, i. e. the condition of the existence of a
generalized value is not invariant with respect to duality, when the players change places.

{\bf Example 2}

Consider the antagonistic game  $G=(X,Y,\Bbb R \sp 2,\omega, F)$ with vectorial payoffs (i.e with respect to
component-wise ordering), where $X=\{x\sb 1 ,x\sb 2 \}, \; Y=\{y\sb 1 ,y\sb 2 \} $ and the function $F$ is given
by the following table:
$$
\vbox{ \halign{&\vrule#&
 \strut\quad\hfil#\quad\cr
&$F$&&$y_1$&&$y_2$&\cr \noalign{\hrule} &$x_1$&&(2,1)&&(1,3)&\cr \noalign{\hrule} &$x_2$&&(4,0)&&(0,5)&\cr } }
\;.$$

We see that
$$V\sb {x\sb 1} =\{(p,q)\in \Bbb R \sp 2 \vert \; (p,q)\leq (2,1)\and (p,q)\leq (1,3)\}=\{(p,q)\in \Bbb R \sp 2
\vert \; (p,q)\leq(1,1)\};$$
$$V\sb {x\sb 2} =\{(p,q)\in \Bbb R \sp 2 \vert \; (p,q)\leq (4,0)\and (p,q)\leq (0,5)\}=\{(p,q)\in \Bbb R \sp 2
\vert \; (p,q)\leq(0,0)\}.$$

Hence $V\sb {x\sb 2}\subseteq V\sb {x\sb 1}$, that is $x\sb 1 \geq \sp {\alpha }x\sb 2$ and therefore $x\sb 1 $
is the $\alpha $-greatest strategy. However, if we "convolve" the vectorial payoffs according the formula
$(p,q)\to p+q $ we receive the following game:
$$
\vbox{ \halign{&\vrule#&
 \strut\quad\hfil#\quad\cr
&$F$&&$y_1$&&$y_2$&\cr \noalign{\hrule} &$x_1$&&3&&4&\cr \noalign{\hrule} &$x_2$&&4&&5&\cr } } \;.$$

In this game,  $x\sb 2$ is the  $\alpha$-greatest strategy but not $x\sb 1$. The new game is an homomorphic
image of the previous one. Thus for the games with ordered outcomes, an $\alpha$-optimal strategy can be changed
by passing on to homomorphic images.

{\bf Example 3}

With the preceding example notation, let the antagonistic game  $G$  have the following function $F$:
$$
\vbox{ \halign{&\vrule#&
 \strut\quad\hfil#\quad\cr
&$F$&&$y_1$&&$y_2$&\cr \noalign{\hrule} &$x_1$&&(2,1)&&(1,2)&\cr \noalign{\hrule} &$x_2$&&(3,1)&&(1,3)&\cr } }
\;.$$

Here the strategy $x\sb 2$ Pareto-dominates the strategy $x\sb 1$ but it holds $V\sb {x\sb 1}=V\sb {x\sb 2}$.
Thus it may be that two strategies are equivalent with respect to the $\alpha$-domination and one of them
strongly Pareto-dominates the other.

{\bf Example 4}

Consider the game  $G=(X,Y,A,\omega, F)$ with ordered outcomes, where the ordered set $(A,\omega)$ is the
following 5-elements lattice: $\{0< a,b,c <1\}$, and with realization function $F$ given by the table:
$$
\vbox{ \halign{&\vrule#&
 \strut\quad\hfil#\quad\cr
&$F$&&$y_1$&&$y_2$&\cr \noalign{\hrule} &$x_1$&&b&&1&\cr \noalign{\hrule} &$x_2$&&c&&a&\cr \noalign{\hrule}
&$x_3$&&1&&0&\cr } } \;.$$

Here the strategy $x\sb 1$ is the $\alpha$-greatest strategy. Extend the set $A$ up to $\bar A$ by adding new
elements $f,g \;\text{and}\; h$ and extend the order on $\bar A$ by setting $\;0 < f,g,h;\quad g<a,c; \quad h<
b,c.\;$ In the obtained game $\bar G $, we have $\; \bar V\sb {x\sb 2} =\{0, g \},\; \bar V\sb {x\sb 1}=\{0,
b\}, \; \bar V\sb {x\sb 3} = \{0\},\quad$ hence the strategies $x\sb 1$ and $x\sb 2$ are incomparable with
respect to the $\alpha$-domination. Thus the condition to be an $\alpha$-greatest strategy is not preserved by
adding some non-realized outcomes.

On the base of shortcomings mentioned above, we can make a conclusion that the set of all outcomes that are
"guaranteed" for the player {\bf 1} by a strategy is not a good characteristic of it. There is more fine
characteristic based on another approach. Let $(A, \omega )$ be an ordered set. It is possible to extend the
relation $\omega $ up to set of all subsets of $A$ by the following two ways:
$$A\sb 1 \leq \sp {\overline{\omega}} A\sb 2 \Leftrightarrow (\forall a\sb 1 \in A\sb 1 )
(\exists a\sb 2 \in A \sb 2 ) a\sb 1 \leq \sp {\omega}a\sb 2 ;$$
$$A\sb 1 \leq \sp {\widetilde{\omega}} A\sb 2 \Leftrightarrow (\forall a\sb 2 \in A\sb 2 )
(\exists a\sb 1 \in A \sb 1 ) a\sb 1 \leq \sp {\omega}a\sb 2 .$$

Applying these relations to the set of the kind $A\sb x =\{F(x,y)\vert \; y\in Y \}$ we obtain notions of so
called $\beta -domination$  and {\it dual} $\beta -domination$. The explicit forms of them are the following
ones:
$$x\sb 1 \leq \sp {\beta } x\sb 2 \Leftrightarrow (\forall y\sb 1 \in Y )
(\exists y\sb 2 \in Y ) F( x\sb 1 ,y\sb 2)\leq \sp {\omega}F(x\sb 2 ,y\sb 1) ;$$
$$x\sb 1 \leq \sp {\bar{\beta} } x\sb 2 \Leftrightarrow (\forall y\sb 1 \in Y )
(\exists y\sb 2 \in Y ) F( x\sb 1 ,y\sb 1)\leq \sp {\omega}F(x\sb 2 ,y\sb 2) ;$$

These relations agree with the approach of Berge (1957) to sort out the players into active and passive ones. An
active player wants to have at least one good outcomes and a passive one does not want to have a bad outcome.

$\beta -domination$ is interrogated in Rozen (2001), it coincides with relation $R\sb {g\sb 2}$ introduced by
Podinivski (1979).
%%%%%%%%%%%%%%%Preference  functor

\head 2 Preference functor
\endhead
The main model under consideration is DMP of the following kind
$$ G= (X,Y,A, \omega ,F),\tag 1$$
where X,Y,A are arbitrary (nonempty) sets, $\omega \subseteq A\times A $ is a (partial) order relation, $F:
X\times Y \to A$ is a map. This model is interpreted in the following way. The set $X$ is a set of alternatives
(strategies) of a decision maker, the set $Y$ is a set of states and $A$ is a set of outcomes (results). The
order relation $\omega $ expresses the preferences of the decision maker. And the map $F$ is a realization
function (although its values are not assumed to be real numbers). We assume the following underlying conditions
(axioms) for any derived preference relation $\rho \subseteq X\sp 2$:

\roster

(A1) the relation $\rho$ is a preorder relation on $X$;

(A2) the relation $\rho$ contains the Pareto-domination relation;

(A3) a "strengthening" of the relation $\omega $ implies a "strengthening" of the relation $\rho$, in other
words, by adding new outcomes or new comparable pairs to $\omega $, the existing preferences of alternatives
have to be saved.

\endroster
We note that the axioms (A1)-(A3) are from the valuable point of view the minimal system of natural
requirements. The mathematical sense of axioms A1 and A2 is clear, and an exact form of the axiom A3 will be
given below on the base of category theory.

Fix sets $X$ and $Y$ and consider the category $\Cal G (X,Y)$ which has as objects the DMPs $ G= (X,Y,A, \omega
,F)$, and a morphism from $G= (X,Y,A, \omega ,F)$ to $G\sp 1= (X,Y,B, \delta ,H)$ is a map $f: A\to B$
satisfying the following two conditions:
$$ a\sb 1\le \sp {\omega} a\sb 2 \Rightarrow f(a\sb 1 ) \le \sp {\delta } f(a\sb 2);\tag 2$$
$$f\circ F= H.\tag 3$$
The condition (2) means that $f$ is an isotonic map of the ordered set $(A,\omega )$ to the ordered set
$(B,\delta )$ and the condition (3) means that $f$ transfers the function $F$ to the function $H$. In other
words, we consider a DMP (1) as a structure with base set $A$, one binary relation $\omega $ and a family of
0-ary operations $\{F(x,y)\}\sb{(x,y)\in X\times Y}$, and the defined morphisms are usual homomorphisms of
structures (see Gr\"{a}tzer (1979)).

Further, a "strengthening" of preferences on the set of outcomes of the given DMP $ G= (X,Y,A, \omega ,F)$ can
be formalized as an isotonic map $f: A\to B$ of the ordered set $(A,\omega ) $ to an ordered set $(B,\delta)$,
and as a result a new DMP $G\sp 1= (X,Y,B, \delta ,H)$ appears, where $ H=f\circ F$. Thus a "strengthening" is a
morphism of the category $\Cal G (X,Y)$.

Denote by $\Cal Q (X)$ the category which objects are preorder relations on $X $, and for
 $\rho \subseteq X\times X $ and  $\sigma \subseteq X\times X $ the set $Hom (\rho ,\sigma )$ contains only
one morphism if  $\rho \subseteq \sigma $ and is empty otherwise. Since our aim is to associate with every DMP a
preorder relation $\rho \subseteq X\times X $ (a derived preference relation) in such a way that the condition
(A3) is satisfied, we have to consider a map $\Psi $ which assigns to every $\Cal G (X,Y)$-object $G$ a preorder
relation $\Psi (G)$ on $X$ in such way that for every morphism $G \to G\sp 1$ the inclusion $\Psi (G) \subseteq
\Psi (G\sp 1) $ is satisfied. Thus we have to construct functors from the category $\Cal G (X,Y)$ to the
category $\Cal Q (X)$. This is the mentioned category-theoretical formalization of the axiom A3.

Consider the set $A\sp Y $ of all maps from $Y$ to $A$. For every order relation $\omega$ on the set $A$, we
have the order relation $\hat {\omega }$ on $A\sp Y $ : $(\alpha ,\beta ) \in \hat {\omega} \Leftrightarrow
(\forall y\in Y )\;\alpha (y) \le \sp {\omega } \beta (y) $. As it is known, we can associate with every map
$F:X\times Y \to A $ the map $F\sp * :X\to A\sp Y \,$, that assigns to every element $x\in X$ the map $F\sb x
\sp * : Y\to A $ by the rule:  $F\sb x \sp * (y)= F(x,y)$  for all $ y\in Y$.

Let $P$ be a map that assigns to every DMP $G= (X,Y,A, \omega ,F)$  the relation $F\sp {*-1} (\hat {\omega }) =
F\sp {*-1} \circ \hat {\omega } \circ F\sp * $. We have

$$P(G) =\{(x\sb 1 ,x\sb 2 )\vert \; (\forall y\in Y)F(x\sb 1,y)\le \sp {\omega } F(x\sb 2 ,y)\}.\tag 4 $$
It is exactly the Pareto-domination for $G$. Let $f: A\to B$ be a morphism from $G= (X,Y,A, \omega ,F)$ to $G\sp
1= (X,Y,B, \delta ,H)$. Then according to definition of morphisms, we have:

$$(x\sb 1 ,x\sb 2 )\in P(G)\Rightarrow (\forall y\in Y)(f\circ F)(x\sb 1,y)\le \sp {\delta } (f\circ F)(x\sb 2
,y)\Rightarrow$$
$$\Rightarrow (\forall y\in Y)\,H(x\sb 1,y)\le \sp {\delta } H(x\sb 2 ,y)\Rightarrow (x\sb 1 ,x\sb 2 )\in P(G\sp
1),$$ that is $P$ is a functor from the category $\Cal G (X,Y)$ to the category $\Cal Q (X)$. We call it the
{\it Pareto functor}.

%%% Definition of a preference functor

\proclaim {Definition 1} A functor $\Psi :\Cal G (X,Y) \to \Cal Q (X)$ is called a preference functor if
$P(G)\subseteq \Psi (G)$ for every DMP $G= (X,Y,A, \omega ,F)$ .
\endproclaim

From the point of view of this definition, the conditions (A1)-(A3) mean that we deal with preference functors
from the category $\Cal G (X,Y)$ to the category $\Cal Q (X)$.

Our definition of a preference functor is very general. All preference relations that are considered in the
literature indeed satisfy more one condition which seems to be natural, but we have never seen it in an explicit
form. Roughly speaking this condition means that we choose the same preferences being in similar situations.
Consider two situations $(x\sb 1,y\sb 1)$ and $(x\sb 2,y\sb 2)$ that lead to outcomes $F(x\sb 1,y\sb 1)$ and
$F(x\sb 2,y\sb 2)$ respectively. We can compare these situations by means of order in $A$, namely the second
situation is preferable than the first one if $F(x\sb 1,y\sb 1)\leq \sp {\omega } F(x\sb 2,y\sb 2)$. This
approach leads to a relation in the set of states for every pair of strategies $(x\sb 1,x\sb 2)$. It seems to be
naturally to take into account this relation. Below we give the exact definition.

\proclaim {Definition 2} The state-preference for the pair $(x\sb 1 ,x\sb 2)$ with respect to DMP $G$ is the
following relation on the set of states $Y$ :
$$\rho \sb G (x\sb 1 ,x\sb 2)=
\{(y\sb 1 ,y\sb 2)\vert \;F(x\sb 1 ,y\sb 1)\leq \sp {\omega }F(x\sb 2 ,y\sb 2)\}.\tag 5$$

\endproclaim

\proclaim {Definition} A preference functor  $\Psi :\Cal G (X,Y) \to \Cal Q (X)$ is called regular if it
satisfies the following condition:

(A4)

if  \quad $\rho \sb G (x\sb 1 ,x\sb 2)= \rho \sb {G'} (x\sb 3 ,x\sb 4) $ then
$$ (x\sb 1 ,x\sb 2)\in \Psi (G) \Leftrightarrow (x\sb 3 ,x\sb 4)\in \Psi (G')\tag 6$$
\endproclaim

Now we show that so called Ordinal Invariance Axiom (OIA) (see {\bf QDT }) is more stronger that A4.  OIA means
in our notations that for $x\sb 1 ,x\sb 2, x\sb 3 ,x\sb 4 \in X$, if $\{y\in Y \vert \;F(x\sb 1 ,y)\leq \sp
{\omega }F(x\sb 2 ,y)\}=\{y\in Y \vert \;F(x\sb 3 ,y)\leq \sp {\delta }F(x\sb 4 ,y)\}$ then $ (x\sb 1 ,x\sb
2)\in \Psi (G) \Leftrightarrow (x\sb 3 ,x\sb 4)\in \Psi (G')$.

Let $\rho \sb G (x\sb 1 ,x\sb 2)= \rho \sb {G'} (x\sb 3 ,x\sb 4) $. The fact $F(x\sb 1 ,y)\leq \sp {\omega
}F(x\sb 2 ,y)$ means that $(y,y)\in \rho \sb G (x\sb 1 ,x\sb 2)$. The last one is equal to $(y,y)\in \rho \sb G
(x\sb 3 ,x\sb 4)$, that is equal for one's turn to $F(x\sb 1 ,y)\leq \sp {\delta }F(x\sb 2 ,y)$. Hence  $\{y\in
Y \vert \;F(x\sb 1 ,y)\leq \sp {\omega }F(x\sb 2 ,y)\}=\{y\in Y \vert \;F(x\sb 3 ,y)\leq \sp {\delta }F(x\sb 4
,y)\}$ and OIA implies that $ (x\sb 1 ,x\sb 2)\in \Psi (G) \Leftrightarrow (x\sb 3 ,x\sb 4)\in \Psi (G')$. Thus
OIA implies A4.
%%%%%%%%%%%%%%%%%%%%%%%%%%%%%%%%%%%%%%%%%%%%%%%

%%%%%%%%%%%%%%%%%%%%%Main result
\head 3. Main results
\endhead
In this section, we give a method of constructing of all regular preference functors $\Psi :\Cal G (X,Y) \to
\Cal Q (X)$ and consider connections between them. It is obvious that the required construction is not to depend
on sets of strategies $X$ and realization functions $F$ because the last ones only choose and rename some maps
from $Y$ to $A$. We should mention that some authors consider from the very beginning preorder relations on
$A\sp Y $, the sets of acts of a decision maker (see for example Barth\'{e}lemy et al. (1982), Dubois et al.
(2002)). But we do not assume in contrast to this point of view that all mappings are possible acts. However by
constructing of preference functors, we consider preorder relations on the set $A\sp Y $ without connections
with functions $F$. The following fact gives a reason for such approach.

Recall that every map $f:A\to B$ determines a map $\tilde f:A\sp Y \to B\sp Y $ under formula: $\tilde f
(\varphi )= f\circ \varphi $  for all  $\varphi : Y\to A $. Denote by $\Cal Q $ the category of all preordered
sets and isotonic maps between them.

\proclaim {Proposition 1}

(1) Let $\Psi :\Cal G (X,Y) \to \Cal Q (X)$  be a preference functor. Consider a map  $\Xi$ that assigns 1) to
every DMP $G= (X,Y,A, \omega ,F)$ the preorder relation on the set $A\sp Y$ generated by $ F\sp * \circ \Psi (G)
\circ F\sp {*-1}$ and  $\hat {\omega }$ and 2) to every morphism of DMPs $f$ the map $\tilde f$.  Then $\Xi$ is
a functor from the category $\Cal G (X,Y)$ to the category $\Cal Q $.

2) Assume that we have a functor $\Xi$ from the category of all ordered sets to the category $\Cal Q $ which
assigns to every ordered set $(A,\omega )$ a preorder relation $\Xi(\omega )$ on the set $A\sp Y$ and such that
the following two conditions are satisfied:

(i) $\hat {\omega }\subseteq \Xi(\omega ) $ ;

(ii) if $f: A\to B $ is an isotonic map of the ordered set $(A,\omega )$ into the ordered set $(B,\delta )$ then
$\Xi (f)= \tilde f :A\sp Y \to B\sp Y $.

Then the map $\Psi$ which assigns to every DMP $G= (X,Y,A, \omega ,F)$ the preorder relation $F\sp {*-1}
(\Xi(\omega ))$ determines a preference functor from category $\Cal G (X,Y)$ to the category $\Cal Q (X)$.

\endproclaim

\demo {Proof} (1) Given a DMP $G= (X,Y,A, \omega ,F)$ , consider the kernel of the map $F\sp * $: $\varepsilon
\sb {F\sp *}=\{(x\sb 1,x\sb 2)\vert \;F\sp *\sb{x\sb 1 }=F\sp * \sb{x\sb 2 }\}$. We see that $\varepsilon \sb
{F\sp *}$ is included in the Pareto-domination relation. Since $P(G) \subseteq \Psi (G)$ and $\varepsilon \sb
{F\sp *} \subseteq P(G)$, the inclusion $\varepsilon \sb {F\sp *} \subseteq \Psi (G)$ holds. Therefore the
preorder relation generated by $ F\sp * \circ \Psi (G) \circ F\sp {*-1}$ and $\hat {\omega }$ is equal to $ \hat
{\omega } \circ F\sp * \circ \Psi (G) \circ F\sp {*-1} \circ \hat {\omega } \cup \hat {\omega }$. Now we prove
that the map $\tilde f :a\sp Y \to B\sp Y $ is isotonic. In fact, we have $ \tilde f \circ F\sp * \circ \Psi
(G)\circ F\sp {*-1} \circ \tilde f \sp{-1}= (f\circ F)\sp * \circ \Psi (G)\circ (f\circ F\sp *)\sp {-1} =H\sp *
\circ \Psi (G)\circ H\sp {*-1}.$

It is well known (and quite obvious) that if $f:A\to A$ is the identity map then $\Xi (f) =\tilde f $ is also
the identity map on $A\sp Y$ and that for $f:A\to B, \; g:B\to C $ the equality $\widetilde {g\circ f}=\tilde
{g} \circ \tilde {f}$. Thus $\Xi $ is a functor.

(2) Let $\Xi$ be a functor satisfying the hypotheses.  Let $G= (X,Y,A, \omega ,F)$ and $G\sp 1= (X,Y,B, \delta
,H)$ and $f:G\to G\sp 1 $ be a morphism. Since $f$ is an isotonic map, the map $\tilde f $ also is isotonic and
therefore $(F\sp * \sb {x\sb 1},F\sp * \sb {x\sb 2})\in \Xi(\omega )\Rightarrow (\tilde f(F\sp * \sb {x\sb
1}),\tilde f (F\sp * \sb {x\sb 2}))\in \Xi(\delta )$. Under formula (3), $\tilde f (F\sp * \sb x) = H\sp *\sb x$
holds for arbitrary $x\in X$. Therefore the following implication is true: $(F\sp * \sb {x\sb 1},F\sp * \sb
{x\sb 2})\in \Xi(\omega ) \Rightarrow (H\sp * \sb {x\sb 1},H\sp * \sb {x\sb 2})\in \Xi(\delta )$. Hence $F\sp
{*-1} (\Xi(\omega )) \subseteq H\sp {*-1} (\Xi(\delta ))$. In other words, we can define that $\Psi (f)$ is an
inclusion $\Psi (G) \subseteq \Psi (G \sp 1)$ . It is enough for conclusion that $\Psi $ is a functor. The
condition (i) implies that $P(G)\subseteq \Psi (G)$, i.e. the functor $\Psi $ is a preference functor. \qquad
\qed \

\enddemo

According to Proposition above, we begin to construct functors $\Xi $ from the category of ordered sets to the
category of preordered sets satisfying the conditions for a regular preference functor, that is

1) it assigns to every ordered set $(A, \omega )$ a preorder relation $\Xi (\omega )$ on the set $A\sp Y$ and to
every isotonic map $f: A\to B $ the map $\tilde{f}:A\sp Y \to B\sp Y$;

2) it includes the Pareto functor: $\hat {\omega }\subseteq \Xi(\omega ) $ for all ordered sets $(A, \omega )$ ;

3) it is regular: if $(A, \omega )$ and $(B,\delta )$ are ordered sets and $\psi \sp {-1} \circ \omega \circ
\varphi = \psi \sp {-1} \sb 1 \circ \omega \circ \varphi \sb 1$ holds for  $\varphi ,\psi :Y\to A$ and $\varphi
\sb 1 , \psi \sb 1 :Y\to B$ then $(\varphi ,\psi )\in \Xi (\omega )\Leftrightarrow (\varphi \sb 1 , \psi \sb 1 )
\in \Xi (\delta )$.

We call such a functor also as a regular preference functor.

We recall that a {\it monoid } is a nonempty set $M$ together with an associative binary operation
(multiplication) and an element 1 called an {\it identity }, such that $m1=1 m =m $ for every $m\in M$. The set
$\frak R (Y)$ of all binary relations on $Y$ is a monoid with respect to composition of binary relations and the
identity $\Delta \sb Y $ (see Introduction). We call a subset $\frak A $ of $\frak R (Y) $ {\it a closed
submonoid } if it satisfies the following conditions: \roster

1) $\alpha ,\; \beta \in \frak A \Rightarrow  \beta \circ \alpha \in \frak A $;

2) $\Delta \sb Y \in \frak A $;

3) $\alpha \in \frak A \and \alpha \subseteq \beta \Rightarrow \beta \in \frak A $.

\endroster
Consider an ordered set $(A,\omega )$ and choose some closed submonoid $\frak A $ of $\frak R (Y)$. Define a
relation $\frak A (\omega )$  on the set $A\sp Y $ setting for arbitrary $\varphi \sb 1 ,\varphi \sb 2 \in A\sp
Y $ :

$$(\varphi \sb 1 ,\varphi \sb 2) \in \frak A (\omega ) \Leftrightarrow
\{(y\sb 1, y\sb 2)\in Y\sp 2 \vert \; \varphi \sb 1 (y\sb 1) \le \sp {\omega } \varphi \sb 2 (y\sb 2)\} \in
\frak A . \tag 7 $$

It can be expressed also by the formula :
$$(\varphi \sb 1 ,\varphi \sb 2) \in \frak A (\omega ) \Leftrightarrow
\varphi \sb 2 \sp {-1} \circ \omega \circ \varphi \sb 1 \in \frak A .\tag 8$$

It is easy to see that $\frak A (\omega )$ is a preorder relation on $A\sp Y $. Indeed, for every $\varphi : Y
\to A$ it holds $\Delta \sb Y \subseteq \varphi \sp {-1} \circ \omega \circ \varphi $ , thus we have $ \varphi
\sp {-1} \circ \omega \circ \varphi \in \frak A $, hence $(\varphi , \varphi )\in \frak A (\omega )$. Further,
if $(\varphi  , \psi )\in \frak A (\omega ) $ and $(\psi , \gamma )\in \frak A (\omega ) $, then $(\varphi  ,
\gamma )\in \frak A (\omega )$, because $\gamma \sp {-1} \circ \omega \circ \psi \circ \psi \sp {-1} \circ
\omega \circ \varphi \subseteq \gamma \sp {-1} \circ \omega \circ \varphi $.

\proclaim {Theorem 1} Let $\frak A $ be a closed submonoid of $\frak R (Y)$. The map $\Psi $ that assigns to
every DMP $G =(X,Y,A,\omega , F)$ the preorder relation $F\sp {*-1} (\frak A (\omega))$ determines a regular
preference functor from the category $\Cal G (X,Y)$ to the category $\Cal Q (X)$.
\endproclaim
\demo {Proof}

Let  $(A, \omega )$ and $(B, \delta )$ be ordered sets. Let  $f: A\to B $ be an isotonic map. Let further
$(\varphi  , \psi )\in \frak A (\omega ) $. The last one means that $ \psi \sp {-1} \circ \omega \circ \varphi
\in \frak A $. Since  $ \psi \sp {-1} \circ \omega \circ \varphi \subseteq \psi \sp {-1} \circ f \sp {-1} \circ
f \circ \omega \circ f\sp {-1} \circ f \circ \varphi \subseteq (f\circ \psi )\sp {-1} \circ \delta \circ (f\circ
\varphi )$, we have $(f\circ \varphi  ,f\circ \psi )\in \frak A (\delta ) $.

The fact above can be expressed in the following words: every closed submonoid $\frak A $ of $\frak R (Y) $
determines a functor $\Xi $ from the category of ordered sets to the category $\Cal Q$ of preordered sets. This
functor assigns to every ordered set $(A,\omega )$ the preordered set $(A\sp Y ,\frak A (\omega ))$ and to every
isotonic map $f$ of $(A, \omega )$ in $(B, \delta )$ the isotonic map  $\tilde f$ of $(A\sp Y ,\frak A (\omega
))$ into $(B\sp Y ,\frak A (\delta ))$ .

Now we apply the part (2) of Proposition.  The condition (ii) is fulfilled. Verify the condition (i), that is
for all $\varphi\sb 1 , \varphi\sb 2 \in A\sp Y $ : $(\forall y\in Y) \varphi \sb 1 (y)\le\sp {\omega }
\varphi\sb 2 (y) \Rightarrow (\varphi\sb 1 , \varphi\sb 2)\in \frak A (\omega ) $ holds.

Indeed, if $(\forall y\in Y) \varphi\sb 1 (y)\le\sp {\omega } \varphi\sb 2 (y)$ then $\Delta \sb Y \subseteq
\varphi \sb 2 \sp {-1} \circ \omega \circ \varphi \sb 1 $ and therefore $\varphi \sb 2 \sp {-1} \circ \omega
\circ \varphi \sb 1 \in \frak A $.

Then the conclusion of Proposition means that the map $\Psi $ which assigns to every DMP $G =(X,Y,A,\omega , F)$
the preorder relation $F\sp {*-1} (\frak A (\omega))$ and to every morphism $f: G= (X,Y,A, \omega ,F) \to G\sp
1= (X,Y,B, \delta ,H)$ the inclusion $F\sp {*-1} (\frak A(\omega ))\subseteq H\sp {*-1} (\frak A(\delta ))$ is a
preference functor from the category $\Cal G (X,Y)$ to the category $\Cal Q (X)$.

Further, we have for two DMPs $ G= (X,Y,A, \omega ,F)$ and $ G\sp 1= (X,Y,B, \delta ,H)$ :
$$(x\sb 1 ,x\sb 2 )\in \Psi (G) \Leftrightarrow (F\sp * \sb{x\sb 1},F\sp * \sb{x\sb 2})\in \frak A (\omega ),$$
$$(x\sb 3 ,x\sb 4 )\in \Psi (G\sp 1) \Leftrightarrow (H\sp * \sb{x\sb 3},H\sp * \sb{x\sb 4})\in \frak A (\delta ).$$
If $\rho \sb G (x\sb 1 ,x\sb 2)= \rho \sb {G\sp 1} (x\sb 3 ,x\sb 4) $, then according to (5) right sides of
these formulas are equivalent and hence $ (x\sb 1 ,x\sb 2)\in \Psi (G) \Leftrightarrow (x\sb 3 ,x\sb 4)\in \Psi
(G\sp 1)$, that is (6) is fulfilled. Thus the functor $\Psi $ is regular.\qquad \qed \

\enddemo

Theorem~1 shows that choosing some closed submonoid $\frak A $ of the monoid of all binary relations on the set
$Y$ we can construct derivative preference relations for every DMP from $\Cal G (X,Y)$. The next theorem shows
that this is a way to obtain all regular preference functors.

\proclaim {Theorem 2}

Let $\Xi $  be a regular preference functor from the category ordered sets to the category $\Cal Q $.  Denote by
$\frak A $ the set of all binary relations on the set $Y$ of the form $\psi \sp {-1 } \circ \omega \circ \varphi
$ where $(\varphi ,\psi )\in \Xi(\omega )$. Then $\frak A $ is a closed monoid of relations and $\frak A (\omega
)= \Xi (\omega )$ for every ordered set $(A,\omega )$.
\endproclaim
\demo {Proof } Consider the trivially ordered set $(Y, \Delta \sb Y )$. Since for every $\varphi :Y\to Y $,
 $(\varphi ,\varphi )\in \Xi ( \Delta \sb Y )$ holds, we have $\varphi \sp {-1 } \circ \Delta \sb Y \circ \varphi
\in \frak A$. For the case $\varphi $ is the identity map, we obtain that $\Delta \sb Y \in \frak A $.

Let $\rho ,\sigma \in \frak A$. Then there are ordered sets $(A,\omega )$ and $(B,\delta )$ and there are maps
$\varphi \sb 1 , \psi\sb 1 :Y\to A $, $\varphi \sb 2 , \psi \sb 2 :Y\to B $ such that $(\varphi \sb 1 ,\psi \sb
1)\in \Xi(\omega )$, $(\varphi \sb 2 ,\psi \sb 2)\in \Xi(\delta )$ and $\rho = \psi \sb 1\sp {-1 } \circ \omega
\circ \varphi \sb 1$, $\sigma = \psi \sb 2\sp {-1 } \circ \omega \circ \varphi \sb 2 $. It can be assumed that
sets $A$ and $B$ have not common elements. Consider the relation $\tau = \omega \cup \delta \cup \delta \circ
\varphi \sb 2 \circ \psi \sb 1 \sp {-1} \circ \omega $
 on the set  $A\cup B $. It is easy to see that $\tau $ is an order relation on this set.
For each $y\in Y $ we have
$$(\psi \sb 1 (y) , \varphi \sb 2 (y)) \in \delta \circ \varphi \sb 2 \circ \psi \sb 1
\sp {-1} \circ \omega \subseteq \tau .$$

Consider two inclusions $f:A \to A\cup B $ and $g:B\to A\cup B $. Under definition $\tau $, they are isotonic.
We conclude from the formula above,  that $(f\circ \psi \sb 1 ,g \circ \varphi \sb 2) \in \Xi (\tau )$. But
since the maps are isotonic, we have also $(f\circ \varphi \sb 1 ,f \circ \psi \sb 1) \in \Xi (\tau )$ and
$(g\circ \varphi \sb 2 ,g \circ \psi \sb 2) \in \Xi (\tau )$. Therefore $(f\circ \varphi \sb 1 ,g \circ \psi \sb
2 ) \in \Xi (\tau )$. Further we have:
$$\psi \sb 2\sp {-1 } \circ g\sp {-1} \circ \tau \circ f \circ \varphi \sb 1 =
\psi \sb 2\sp {-1 } \circ \delta \circ \varphi \sb 2 \circ \psi \sb 1 \sp {-1} \circ \omega \circ \varphi \sb 1
= \sigma \circ \rho $$ and hence $\sigma \circ \rho \in \frak A $. Thus we have proved that $\frak A $ is a
monoid.

Now we have to proof that this monoid is a closed submonoid of $\frak R (Y)$. Firstly observe a trivial fact
that every relation $\sigma \subseteq Y\times Y $ can be represented in the form: $\sigma =\psi \sp {-1 }\circ
 \omega \circ \varphi $ for some ordered set $(A, \omega )$ and maps $\varphi ,\psi :Y \to A $. Indeed, let $Y\sb
 1$ and $Y\sb 2$ be two copies of $Y$ without common elements. Let $\varphi :Y\to Y\sb 1$ and $\psi :Y\to Y\sb 2$
be the  corresponding identification maps. Define $\omega = \Delta \sb {Y\sb 1} \cup \Delta \sb {Y\sb 2}\cup
\psi \circ  \sigma \circ \varphi \sp {-1}$. It is true that $\omega $ is an order relation on $A=Y\sb 1 \cup
Y\sb 2 $ because a pair of different elements belongs to $\omega $ only if the first of them belongs to $Y\sb 1
$ and the second one belongs to $Y\sb 2$. It is easy to see that $\sigma =\psi \sp {-1 } \circ \omega \circ
\varphi $.

Let $\rho \subseteq \sigma $ and $\rho \in \frak A $. The last one means that there exist $\varphi \sb 1, \psi
\sb 1 :Y\to B $ and order relation $\delta $ in $B$ such that $\rho = \psi \sb 1\sp {-1 } \circ \delta \circ
\varphi \sb 1 $ and $(\varphi \sb 1, \psi \sb 1) \in \Xi (\delta )$. As it was mentioned above,  $\sigma = \psi
\sp {-1 } \circ \omega \circ \varphi$ for some maps $\varphi , \psi  :Y\to B $ and order relation $\omega $ in
$B$. Consider the ordered set $(A\times B , \omega \times \delta )$ where $((a\sb 1 ,b\sb 1 ), (a\sb 2 ,b \sb 2
)) \in \omega \times \delta \Leftrightarrow (a\sb 1 ,a \sb 2 )\in \omega \and (b\sb 1 ,b \sb 2 ) \in \delta $.
We have two natural maps $\alpha ,\beta : Y \to A\times B$ where $\alpha (y) = (\varphi (y),\varphi \sb 1 (y))$
and $\beta (y) = (\psi (y),\psi \sb 1 (y))$.  It is easy to see that $\beta \sp {-1 } \circ (\omega \times
\delta ) \circ \alpha = \rho $. Since the functor $\Xi$ is regular, the last equality means that $(\alpha ,\beta
) \in \Xi (\omega \times \delta )$. Using the isotonic map $\pi :A\times B \to A $, projection on $A$, we obtain
that $(\varphi ,\psi )\in \Xi (\omega )$, and hence $\sigma \in \frak A $. Therefore the monoid $\frak A $ is
closed.

Now let $(\varphi ,\psi )\in \frak A (\omega )$ for some ordered set $(A,\omega )$. It is equivalent according
to (8) that $ \rho=\psi \sp {-1 } \circ \omega \circ \varphi \in \frak A$. The last fact means that there exist
$\varphi \sb 1, \psi \sb 1 :Y\to B $ and order relation $\delta $ in $B$ such that $\rho = \psi \sb 1\sp {-1 }
\circ \delta \circ \varphi \sb 1 $ and $(\varphi \sb 1, \psi \sb 1) \in \Xi (\delta )$. Since the functor $\Xi $
is regular, it is equivalent to $(\varphi ,\psi )\in \Xi (\omega )$. Therefore $\frak A (\omega )= \Xi (\omega
)$. \qquad \qed \

\enddemo

%%%%%section 4 . Exaples and applications
%%%%%%%%%%%%%%%%%%%%%%%%%%%%%%%%%%%%%%%%%%%%%%%%%%??????????????????????????????
\head 3. Examples and applications.
\endhead
\bigskip
Thus to observe all closed preference functors from the category $\Cal G (X,Y) $, we have to observe all closed
monoids of relations on the set $Y$. In reality, it is sufficient to choose a submonoid of $\frak R (Y)$ and
then consider all relations which include its members. All closed submonoids of $\frak A $ form a lattice with
respect to inclusion. This lattice is complete. The least element of this lattice is submonoid of all reflexive
relations and the greatest one is $\frak R (Y)$ itself. We show below that the submonoid of all reflexive
relations gives the Pareto-domination. Clearly, $\frak R (Y)$ gives the greatest preference relation, i. e. the
complete relation on $X$: $X\times X$. If $\frak A $ is a closed submonoid of $\frak R (Y)$, then the set of all
relations $\varrho \sp {-1}$ for $\varrho \in \frak A $ also is a submonoid of $\frak R (Y)$. Thus to every
preference functor, there is the dual one. If they coincide, we have a self-dual preference functor. Below we
give some useful examples.

{\bf Examples}

1. Consider the set of all reflexive relations on the set $Y$. Clearly, it is a submonoid of $\frak R $. Let
$\frak A$ be this submonoid. The corresponding derivative preference relation is the Pareto-domination. Indeed,
$$(\alpha ,\beta )\in \frak A (\omega)\Leftrightarrow \Delta \sb Y \subseteq
\beta \sp {-1} \circ \omega \circ \alpha \Leftrightarrow (\forall y\in Y )\alpha (y)\le \sp {\omega} \beta (y)
\Leftrightarrow (\alpha ,\beta )\in \hat {\omega }.\tag 9$$

2. We call a relation on $Y$  surjective if its second projection is equal to Y. Clearly, all such relations
form a submonoid of $\frak R $. Let $\frak A$ be the submonoid of all surjective relations. The corresponding
derivative preference relation is the $\beta $-domination defined in Section 1. Indeed,
$$(\alpha ,\beta )\in \frak A (\omega)\Leftrightarrow
pr\sb 2 ( \beta \sp {-1} \circ \omega \circ \alpha)=Y \Leftrightarrow (\forall y\in Y )(\exists y'\in Y)\;
\alpha (y')\le \sp {\omega} \beta (y).\tag 10$$ Thus
$$ (x\sb 1, x\sb 2)\in F\sp {*-1} (\frak A (\omega))\Leftrightarrow (\forall y)(\exists y')
F( x\sb 1,y') \le \sp{\omega} F(x\sb 2,y ).\tag 11$$

3. Dually to the previous example, let $\frak A$ be the submonoid of all everywhere defined relations (it means
the first projection of a relation is equal to Y). The corresponding derivative preference relation is the dual
$\beta $-domination. The proof is the same that above.

4. Recall that a {\it filter } on the set $Y$ is a set $\frak F $ of non-empty subsets of $Y$ satisfying the
following conditions: 1)$Y\in \frak F$, 2) $A,B \in \frak F $ implies $A\cap B \in \frak F$ and 3) $A\in \frak
F$ and $A\subseteq B$ implies $B\in \frak F$.

Let $\frak F $ be a filter on the set $Y$. Define $\frak A =\{\lambda \subseteq Y\times Y \vert \; pr\sb {\Delta
}\in \frak F \}$. The corresponding derivative preference relation is the preference according to the filter. It
can be expressed as follows: a strategy $x\sb 2 $ is more preferable than a strategy $x\sb 1 $ if the set of all
$y$ for which $F( x\sb 1,y ) \le \sp{\omega} F(x\sb 2,y )$ belongs to the filter $\frak F $. The well known
interpretation: "the majority vote for".

A filter $\frak F $ is called {\it principal } if it is of the following form: $\frak F=\{P\subseteq Y \vert \;
P\sb 0 \subseteq P $, where $P\sb 0 $ is a fixed subset of $Y$. One can consider the set $P\sb 0 $ as a system
of dominators. If particularly $P\sb 0 $ consists of one state only, this state can be considered as a indicator
(or a dictator): a strategy $x\sb 1$ is preferable than $x\sb 2$ if and only if the strategy $x\sb 1$ gives for
this state the result better than $x\sb 2$. If the set $Y$ is finite, every filter is principal.

5. Let $\sigma $ be an idempotent relation on the set $Y$, that is $\sigma \circ \sigma =\sigma $. It means that
a decision maker has some special preferences on the set of states. Denote by $\frak A$ the closed monoid
generated by $\sigma $. It consists of reflexive binary relations and all binary relations containing $\sigma $.
We have
$$ (\alpha ,\beta )\in \frak A (\omega)\Leftrightarrow (\alpha ,\beta )\in \hat {\omega} \lor
(\forall y\sb 1, y\sb 2\in Y)((y\sb 1,y\sb 2)\in \sigma \Rightarrow  \alpha (y\sb 1)\le \sp {\omega} \beta (y\sb
2)).\tag 12$$

According to this we obtain a new derived preference relation.
\bigskip

We recall that Example 3 given in Section 1 presents two strategies equivalent with respect to
$\alpha$-domination such that one of them strongly Pareto-dominates the other. Although $\alpha$-domination can
not be realized by means of a preference functor, such undesirable event can appear, for example if one chooses
the greatest monoid $\frak A $.

\proclaim {Definition 3} Let $G$ be a DMP. Let $P(G)\sb {str}$ denote the strong Pareto-domination for $G$. A
preference functor $\Psi :\Cal G (X,Y) \to \Cal Q (X)$ is called suitable for the  $G$  if the following
condition is satisfied:

(A5)
$$ \Psi (G) \cap P(G)\sb {str} \sp {-1} = \emptyset.$$
\endproclaim

The condition above means that it is impossible that two strategies $x\sb 1$ and $x\sb 2$ are equivalent with
respect to preference relation $ \Psi (G)$ but $x\sb 2$  strong Pareto-dominates  $x\sb 1$ . It is obvious that
Pareto preference functor is suitable for all DMPs. In general case a preference functor can be suitable for
some DMPs and non-suitable for another ones.

\proclaim {Proposition 2} Let $\frak A $ be a closed submonoid of $\frak R (Y)$. If every relation $\rho \in
\frak A $ has a fixed point $y$, i. e.  $(y,y)\in \rho $, then the preference functor $\Psi$ determined by
$\frak A $ is suitable for every DMP.

\endproclaim
\demo {Proof} Let $G=(X,Y,A,\omega,F)$ be a DMP. Suppose that $(\alpha ,\beta )\in \frak A (\omega)$. It means
that $  \beta \sp {-1} \circ \omega \circ \alpha \in \frak A $. Under hypotheses, this relation has a fixed
point $y\sb 0$. For this point $\alpha (y\sb 0) \leqq \sp {\omega } \beta (y\sb 0)$ holds. Hence it is
impossible that $(\forall y )\beta  (y)< \sp {\omega } \alpha (y)$. It means that $ \Psi (G) \cap P(G)\sb {str}
\sp {-1} = \emptyset.$ \qquad \qed \

\enddemo

\proclaim {Corollary 1} The preference functor according to a filter on the set $Y$ (example 4) is suitable for
every DMP .
\endproclaim

%Recall that a {\it chain } in an ordered set $(A,\omega)$ is a sequence of elements.

\proclaim {Proposition 3} Let $\frak A $ be a non-universal closed submonoid of $\frak R (Y)$. If
 $ \frak A (\omega) \cap \omega \sb {str} \sp {-1} \not = \emptyset$, then $(A,\omega)$ contains
a chain $a\sb1 <a\sb 2 < ... <a\sb k$ for every positive integer k.
\endproclaim
\demo {Proof} Suppose that there are  $\varphi, \psi \in A\sp Y$ such that  $(\varphi , \psi)\in \frak A (\omega
)$ but $\psi (y) < \sp {\omega } \varphi (y)$ for all $y\in Y$.  Then the relation $\rho =\{(y\sb 1 ,y \sb 2
)\vert \; \varphi (y\sb 1)\leqq \sp {\omega } \psi ( y \sb 2)\} = \psi \sp {-1} \circ \omega \circ \varphi $
belongs to $\frak A $. Hence $\rho \sp k \in \frak A$ for arbitrary natural number $k$. Since $ \varphi \circ
\psi \sp {-1} \subseteq \omega \sb {str}$,  $\rho \sp 2 = \psi \sp {-1} \circ \omega \circ \varphi  \circ \psi
\sp {-1} \circ \omega \circ \varphi  \subseteq \psi \sp {-1} \circ \omega \sb {str} \circ \varphi $, and further
$\rho \sp 3 \subseteq \psi \sp {-1} \circ \omega \sp 2\sb {str} \circ \varphi $, and so long, ..., $\rho \sp k
\subseteq \psi \sp {-1} \circ \omega \sp {k-1}\sb {str} \circ \varphi $. Since submonoid $\frak A $ is closed
$\psi \sp {-1} \circ \omega \sp k\sb {str} \circ \varphi \in \frak A $ for all $k$. Under hypotheses, $\frak A $
is not universal and therefore all these relations are not empty. It gives a sequence $a\sb 1 , ...,a \sb k$ for
every integer $k>0 $ such that $a\sb1 <a\sb 2 < ... <a\sb k$ . \qquad \qed \
\enddemo

\proclaim {Corollary 2} If lengths of all strong chains in the set of outcomes of a DMP $G$ are bounded above,
then every non trivial derivative preference relation is suitable for $G$, particularly for all DMPs with finite
set of outcomes.
\endproclaim

Proposition 3 generalizes the following result obtained in [Rozen (2001), Theorem 2].

\proclaim {Corollary 3} If the ordered set of outcomes for a DMP $G$ satisfies descending (increasing) chain
condition then $\beta$-domination (inverse $\beta $-domination) preference functor is suitable for $G$.
\endproclaim
\demo {Proof} The submonoid $\frak A $ corresponding to $\beta$-domination preference functor consists of all
surjective relations. Suppose that $ \frak A (\omega) \cap \omega \sb {str} \sp {-1} \not = \emptyset$. Then
with the notation of the proof of Proposition 3, the relation $\psi \sp {-1} \circ \omega \sp k \sb {str} \circ
\varphi $  is surjective. It means that for every $y\in Y$ there is $y\sb 1 \in Y$ such that $(\varphi (y\sb 1),
\varphi (y))\in \omega \sb {str}\sp k$. It leads to an infinite descending chain in A. This contradiction shows
that $\beta$-domination preference functor is suitable for $G$.

\enddemo

\head 5. Conclusion
\endhead

The approach suggested in the paper and based on category theory leads to a strong definition of derived
preferences of a decision maker and the obtained results give a method to observe all such preferences and
compare them with each other. The most important point is that we start from minimal conditions and therefore
include all reasonable other approaches but exclude all anomalies mentioned in Section 2. The second point is
that we reduce the problem formulated in Introduction to pure algebraical one, namely to study closed monoids of
binary relations. Indeed, using Theorems 1 and 2 it is possible to build (of course, if the set $Y$ of states is
finite) all regular preference functors for given set $X$ of strategies. One must mention that these theorems
establish connections between regular functors (from the category of all ordered sets to the category of
preordered sets) and closed monoids and are new.

It was said that all closed monoids of relations on the set $Y$ form a lattice. This lattice is a complete
lattice. If a monoid $\frak A\sb 1$ is a submonoid of a monoid $\frak A\sb 2$ and  $\xi \sb 1 $ and $\xi \sb 2 $
are the corresponding derived preferences relations on the set $X$  for a DMP $G$,  then we have $\xi \sb 1
\subseteq \xi \sb 2 $. It can be interpreted that preference $\xi \sb 2 $ is stronger than $\xi \sb 1 $: if
$x\sb 2$ is preferable than $x\sb 1$ in the sense of $\xi \sb 1 $ it is preferable than $x\sb 1$ in the sense of
$\xi \sb 2 $.

In the simple cases, when the set $Y$ is not too large, one can describe explicitly all derivative preference
relations. For example, when the set $Y$ has only two elements there are only 16 binary relations on $Y$, the
derivative preference relations are: two extremal relations (the universal one and the Pareto-domination), two
relations associated with two principal filters, the $\beta$-domination, the inverse $\beta$-domination (four
maximal relations only) and their intersections. In the case when $Y$ contains more elements it is more
complicated computational problem, but is seems that the most important thing is to give a qualitative
description which does not depend on a set of states.

The Pareto-domination is the weakest derivative preference relation. The universal (all strategies are
equivalent) relation is the strongest derivative preference relation. It is very easy to proof that the monoid
of all surjective relations and the monoid of all everywhere defined relations are maximal closed submonoids
(dual atoms). It means that $\beta$-domination and inverse $\beta$-domination are maximal preferences, i.e.
there are no preference relations between each of them and the universal relation.

It may be that a decision maker has some preferences also in the set $Y$. One can take such preferences in
consideration choosing a filter on the set $Y$ or some idempotent relation on it. This leads to the derived
preference relations described in the examples 4 and 5 in the Section 4. The preference relations associated
with principal filters (the case of dictators) are also maximal ones.

Taking intersections of these known monoids we obtain new preference relations that are weaker than their
parents. For example, considering the intersection of $\beta$- and inverse $\beta$-domination (i. e. the monoid
of all everywhere defined surjective relations) gives a new preference relation that can be expressed by the
formula: $x\sb 2$ is preferable than $x\sb 1 $ iff
$$(\forall y\sb 1 )(\exists y\sb 2) F(x\sb 1, y\sb 1)\leq F(x\sb 2, y\sb 2)
\wedge (\forall y\sb 2 )(\exists y\sb 1) F(x\sb 1, y\sb 1)\leq F(x\sb 2, y\sb 2 ).\tag 13$$

On the other hand one can consider minimal preference relations (atoms). Let $y\sb 0 \in Y$ be a fixed element.
Build the relation $\varrho =Y\times Y \setminus \{(y\sb 0 , y\sb 0)\}$. It is a maximal relation and $\varrho
\circ \varrho =Y\times Y $. Thus all reflexive relations and $\varrho $ form a closed monoid that clearly is a
minimal element in the lattice of all closed monoids. This is a way to build all minimal derived preference
relations. The sense of them is obvious: $x\sb 2$ is preferable than $x\sb 1 $ iff  $x\sb 2$ Pareto-dominates
$x\sb 1 $ or $F(x\sb 1, y\sb 1)\leq F(x\sb 2, y\sb 2)$  for all pairs of states $(y\sb 1 , y\sb 2)$ different
from $( y\sb 0 , y\sb 0)$ . It means that if a decision maker wants to ignore for some reason the correlation
between $F(x\sb 1, y\sb 0)$ and $ F(x\sb 2, y\sb 0)$  he has to check the inequality above for all pairs of
states, with the exception of $( y\sb 0 , y\sb 0)$.

We see that there are many different preference functors. But it is not a shortcoming, it is the nature of
things. Practically, we choose one of these functors following supplementary data. For example, one can consider
some structures on the set of states like mentioned above or a probability distribution and use corresponding
preference functors.

\newpage
\head References
\endhead
Barth\'{e}lemy, Cl. Flament, and B. Monjardet.,1982. {\it Ordered sets and social sciences}. In I.Rankin,
editor, {\it Ordered Sets}, 721-787. D. Reidel, Dordrecht.

Berge, C.,1957. {\it Th$\acute e$orie g\'{e}n\'{e}rale des jeux \`{a} n-personnes}, M\'{e}mor. Sci. Math.,no
138. Gauthier-Villars, Paris.

Dubois, D., Fargier, H., Perny P., and Prade H., 2002. {\it Qualitative decision theory: From Savage's axioms to
nonmonotonic reasoning}. Journal of ACM, 49, No 4, 455-495.

Gr\"{a}tzer, G.,1979. {\it Universal Algebra}, 2nd ed., Springer, New York.

Jentzsch, G.,1964. {\it Some thoughts on the theory of cooperative games}. Advances in Game Theory 52, 407-442.

Mac Lane, S.,1971. {\it Categories for the Working Mathematician}. Springer, New York.

Moulin, H.,1981. {\it Th$\acute e$orie des jeux pour l'$\acute{e}$conomie et la politique }(Paris).

Podinovski, V.,1979. {\it Principe of guaranteed results for partial preference relations}. Zh. Vychisl. Mat. i
Mat. Fiz. 19, No 6, 1436-1450 (Russian).

Podinovski, V.,1981. {\it Generalized antagonistic games}. Zh. Vychisl. Mat. mat. Fiz. 21, No 5, 1140-1153
(Russian). English translation in U.S.S.R. Comput. Maths. Math. Phys. Vol.21, No.5, pp.65-79, 1981.

Rozen, V.,1978. {\it $P$-extension of games with quasiordered outcomes}. Mathematical models of behavior
(Russian), pp. 100--115, 126, Saratov. Gos. Univ., Saratov. MR0566719 (81f:90124).

Rozen, V.,2001. {\it Order invariants and the "environment" problem for games with ordered outcomes}. Kibernet.
Sistem. Anal., No 2 , 145-159 (Russian). English translation in Cybernet. Systems Anal. 37 (2001), no. 2,
260--270.

\enddocument